\documentclass[11pt,russian,a4paper,twoside]{article}  %MikTeX
\usepackage[cp1251]{inputenc}                            %Win
\usepackage{amsmath,amsfonts,amssymb,amscd,euscript}
\usepackage[russian]{babel}

\makeatletter 

\def\titlerus{\thispagestyle{empty} { } \vspace{-5mm} \noindent
\raisebox{-37pt}[\headheight][0pt]{\vbox{ \hbox to \textwidth{\hfil
\scriptsize  ВЕСТНИК \; УДМУРТСКОГО \; УНИВЕРСИТЕТА\hfil }
\vspace{2pt} \hrule \vspace{8pt} \hbox to \textwidth{\series \hfil  \issue}
\vspace{30pt} \hbox{УДК \UDK} }} \vspace{ 30pt plus 6pt }} 

\def\titleeng{\vspace{3ex} \hfill Поступила  в редакцию \datereceive \par \vspace{5ex} \par
\noindent \parbox{166mm%130mm
}{\small {\textbf {\textit {\autorseng}}} \par {\bf \articleseng} \par
\vspace{10pt} \par \annotationeng \par \vspace{7pt} \par {\it Keywords}: \keywordseng
\par \vspace{7pt} \par \noindent \small {Mathematical Subject Classifications}: \MSC }
\par \vspace{30pt} \par \small \noindent \contactinformation} 

\def\annotationandkeywordsrus{\noindent {\small \annotationrus \par } \vspace{8pt}
\noindent {\small {\it Ключевые слова}: \keywordsrus} \par \vspace{10pt}}

\@addtoreset{equation}{section}
\renewcommand{\section}{\@startsection{section}{1}{0pt}{1.3ex
plus 1ex minus .1ex}{1.3ex plus .1ex}{\bf\,\S\,}}

\renewcommand{\@begintheorem}[2]{\begin{trivlist}
\item[\hspace{\labelsep}{\bf \mbox{~~~}#1\ #2.}]}
\renewcommand{\@opargbegintheorem}[3]{\begin{trivlist}
\item[\hspace{\labelsep}{\bf \mbox{~~~}#1\ #2 {\rm (#3).}}]}
\renewcommand{\@endtheorem}{\end{trivlist}}

\newtheorem{teo}{Теорема}

\newtheorem{utv}{Утверждение}
\newtheorem{sle}{Следствие}

\newtheorem{zam}{Замечание}
\newtheorem{pr}{Пример}

\newcommand{\doc}{\mbox{Д о к а з а т е л ь с т в о}}
\renewcommand{\@evenfoot}{}
\renewcommand{\@oddfoot}{}

\renewcommand*{\@biblabel}[1]{#1.\hfill}

\newcommand*{\CSep}{.\ }
\renewcommand{\@makecaption}[2]{%
  \vskip\abovecaptionskip
  \sbox\@tempboxa{{\bf #1\CSep}{#2}}%
  \ifdim \wd\@tempboxa >\hsize
  \begin{center}%
    {\footnotesize{\bf #1\CSep}{#2\par}}%
  \end{center}%
  \else
    \global \@minipagefalse
    \hb@xt@\hsize{\hfil\box\@tempboxa\hfil}%
  \fi
  \vskip\belowcaptionskip%
}

\renewcommand{\@evenhead}{\raisebox{0pt}[\headheight][0pt]{\vbox{\hbox to\textwidth{\thepage \strut \hfil
\text{\autorsrus} \hfil } \hrule \vspace{8pt} \hbox to \textwidth{\series \hfil  \issue}}}}

\renewcommand{\@oddhead}{\raisebox{0pt}[\headheight][0pt]{\vbox{\hbox to\textwidth{  \strut \hfil
\text{\articleshortname} \hfil \thepage} \hrule \vspace{8pt} \hbox to \textwidth{\series \hfil  \issue}}}}

\newcommand{\series}{МАТЕМАТИКА}

\newcommand{\issue}{2011. Вып.\,1} %%% Год и номер выпуска.

\newcommand{\autorsrus}{Д.\,В.~Хлопин} %%% на русском языке (точку в конце не ставим)
\newcommand{\autorseng}{D.\,V.~Khlopin} %%% на английском языке (точку в конце не ставим)

\newcommand{\articleshortname}{О расширении игры сближения-уклонения
 на бесконечном промежутке времени}
\newcommand{\articleseng}{On extension   for infinite horizon game of pursuit-evasion}

\newcommand{\UDK}{517.977.57+517.978.4} %%% Проставляет автор!!!
\newcommand{\MSC}{49N70 49N75 49K40} %%% Проставляет автор!!!
\newcommand{\annotationrus}{
 В работе строится расширение конфликтно-управляемых задач на
 бесконечном промежутке. Соответствующее расширение является
 проективным пределом сужений исходной игры на ограниченные промежутки времени.
 Существование максимина в такой расширенной игре эквивалентно
 нечувствительности исходной игры к расширению целевого множества.
 Особое внимание в работе уделяется игре сближения-уклонения
 в паре <<смешанное управление / обобщенное управление>>.
}
\newcommand{\annotationeng}{
  The  extension of a conflict control problem with  infinite  horizon
  is constructed. This extension is
   the projective limit of restricted games.
   Relations between
"sensitivity to target set" and the existence of the optimal control
are studied. Special attention is paid to the pursuit-evasion game
with "joint control\ /\ relaxed control".
}

\newcommand{\keywordsrus}{
расширение игровой задачи, задача сближения-уклонения на бесконечном
промежутке времени, обратные
  спектры.%, %компактно-открытая топология, дискриминация противника.
  }
\newcommand{\keywordseng}{ an extension
for infinity horizon game, %a pursuit-evasion game,
existence of maximin,
  sensitivity of value.
}

\newcommand{\datereceive}{25.10.10} %%% формат дд.мм.гг

\newcommand{\contactinformation}{Хлопин Дмитрий Валерьевич,
к.\,ф.-м.\,н., Институт математики и механики УрО РАН, 620219,
Россия, г. Екатеринбург, ул. Софьи Ковалевской,\,16, E-mail:
khlopin@imm.uran.ru}
\newcommand{\rav}{\stackrel{\triangle}{=}}
\newcommand{\td}[1]{\widetilde{#1}}
\newcommand{\rref}[1]{$(\ref{#1})$}
\newcommand{\mm}[1]{{\bf{#1}}}

\newcommand{\ph}{\varphi}
\newcommand{\pust}{\varnothing}
\newcommand{\bo}{{\hfill {$\Box$}}}
\begin{document}
\sloppy

\titlerus

\begin{flushleft}
{\bf \copyright { \textit { \autorsrus}} \\[2ex]
{ О РАСШИРЕНИИ КОНФЛИКТНО-УПРАВЛЯЕМЫХ ЗАДАЧ \\ НА БЕСКОНЕЧНОМ ПРОМЕЖУТКЕ}%%% ПОЛНОЕ! НАЗВАНИЕ СТАТЬИ
%%% в следующих трех строках Вы проставляете "свои" гранты. Если Ваша работа
%%% не поддержана грантами, эти строки надо закомментировать (лучше стереть)
\footnote{Работа поддержана программой президиума РАН
<<Математическая теория управления>> }}
\end{flushleft}

\annotationandkeywordsrus

%\centerline{\bf Хлопин Д.В.}
%\medskip
%\centerline{\it Институт математики и механики УрО РАН, Россия}
%\centerline{\it khlopin@imm.uran.ru}
%\bigskip
\sloppy
%\begin{abstract}
% В работе исследуются  вопросы существования оптимальных решений в
% игровых задачах на бесконечном промежутке времени. Строятся в
% терминах
% обратных спектров такие расширения исходных вспомогательных задач,
% в которых максимин уже достигается. Для задачи сближения-уклонения
% на бесконечном промежутке времени такое расширение позволяет
% сводить вопрос корректности задачи к вопросу существования
% максимина.
%\end{abstract}
%
% Ключевые слова:
%  расширение игровой задачи, задача сближения-уклонения на бесконечном промежутке времени,
%  существование максимина, устойчивость оптимального
%  результата, топология Виеториса, обратные
%  спектры, компактно-открытая топология, дискриминация противника,
%  extension  for infinity horizon game, existence of optimal control,   compact-open topology.
%a pursuit-evasion game, the projective limit
%%
%
% {\bf Введение.}
\noindent{\bf{Введение.}}

   В теории управления многие исходные постановки задач являются
   некорректными. Самое простое проявление этого ---
   отсутствие оптимального решения, но это может выразиться и в
   скачкообразной зависимости
   оптимального результата от внешних параметров. В случае, когда решение
   отсутствует,
   надо воспользоваться высказанным еще Гильбертом принципом ---
   всего лишь обобщить надлежащим образом понятие решения.
   Соответствующий подход --- расширение исходной задачи ---
   разрабатывался многими авторами, отметим лишь несколько фамилий: для вариационного исчисления
   Л.~Янг  \cite{Yang}, для общих задач управления Р.В.Гамкрелидзе \cite{ga}
   и Дж.~Варга \cite{va}.
%   для задач выпуклого программирования Гольдштейн и Даффин.
%   (ССЫЛКИ)

   Необходимость таких расширений естественно возникает и в дифференциальных играх.
   При доказательстве фундаментальной теорем об альтернативе \cite{ks}
   Н.\,Н.~Красовского и А.\,И.~Субботина ключевым является
   свойство стабильности, предложенное Н.Н.Красовским и
   предусматривающее механизм обобщенной реакции одного из игроков.
   Однако расширения в таких играх не всегда сводятся лишь к введению обобщенных
   программных управлений каждого из игроков, приходится расширять и
   их совместные управления. Например, в одном из подходов,
   предложенных Н.\,Н.~Красовским \cite{ks} к исследованию дифференциальных игр,
   исходная игра подменяется вспомогательной игровой задачей, в
   которой один из игроков заранее знает выбор соперника,
   а следовательно имеет возможность
    рассчитать пучок всевозможных совместимых с таким выбором соперника
    траекторий и выбрать из них любую. В таких игровых задачах с
    дискриминацией соперника фактически игрок выбирает не среди своих
    программных управлений, а среди траекторий как реализаций
    <<совместных>>   управлений. Соответствующие пучки траекторий могут оказываться
    незамкнутыми уже в самых простых игровых задачах, а
    следовательно требовать таких расширений, само множество таких пучков также
    может оказаться незамкнутым.
    Некоторые формализации дифференциальных игр
    сразу приводят к такой дискриминации, например игра в паре
       (чистая стратегия, контрстратегия); см. в этой связи \cite{ks,ou1,iz, 1985}.

    Данная статья идеологически продолжает \cite{kr,2},
     также как и \cite{ChKh} посвящена рассмотрению асимптотических эффектов, возникающих
     во вспомогательных
     игровых задачах с дискриминацией. Общая схема расширения,
     предложенная в этой статье --- уточнение (для задач на
     бесконечном промежутке)  предложенной в \cite{ChKh} схемы расширения  вспомогательных
     игровых задач с дискриминацией.

     Фактически, для задач  на
     бесконечном промежутке сначала исходная игра рассматривается на возрастающей
     последовательности конечных промежутков времени, для каждой такой срезанной игры
     строится на основе метрики Хаусдорфа (фактически по процедуре \cite{ChKh}) её
     расширение (автоматически обладающее  хорошими свойствами),
     после чего берется обратный (проективный) предел этой последовательности расширений.
     То, что при этом получается, может быть создано напрямую   из исходной игры
     при помощи компактно-открытой топологии и топологии Виеториса.
     Для функционалов, непрерывно зависящих от всей траектории,
     максимин так расширенной игры достигается, однако в самом
     интересном случае --- игра сближения-уклонения
     --- соответствующий функционал перестает быть непрерывным
     (в естественной для расширения компактно-открытой топологии), а
     следовательно  существование обобщенных оптимальных
     управлений уже не гарантируется. С другой стороны оказывается, что
     существование таких управлений (в игре-расширении)
     эквивалентно <<грубости>> (устойчивости к малым
     возмущениям целевого множества)
      исходной игровой задачи.

     Полученные результаты иллюстрируются на примере вспомогательной
     игровой задачи сближения-уклонения при формализации
  <<программное смешанное  управление\ /\ программное обобщенное
  управление>>, для которой расширение фактически совпадает с исходной
  игрой.

\noindent{{\bf Определения и обозначения. }}

 Определим $\mm{T}\rav \{t\in\mm{R}\,|\,t\geqslant 0\}.$
 Пусть задано %сепарабельное
  метрическое пространство $\mm{X}.$
   Введем на $C(\mm{T},\mm{X})$ компактно-открытую топологию
   \cite{en};
 в силу \cite[V.3.2]{phillvv} тогда оно метризуемо.
% и некоторое целевое замкнутое множество $M\subset \mm{T}\times\mm{X}.$
 Пусть также задано некоторое множество $\mm{K}\subset C(\mm{T},\mm{X}),$
 причем для всех $t\in \mm{T}$ множества
 $\mm{K}_t\rav \{x|_{[0,t]}\,|\,x\in \mm{K} \}$ компактны.
 Через $exp\,\mm{K}$ обозначим семейство всевозможных непустых
 подмножеств $A$ множества $\mm{K}$ со свойством: для всех $t\in \mm{T}$
 $A_t\rav \{x|_{[0,t]}\,|\,x\in A \}$ замкнуто.
 Пусть  дано некоторое  множество $\mm{V}\subset exp\,\mm{K}.$

%{\bf Функция платы.} Пусть задана некоторая непрерывная функция
% качества $m:\mm{T}\times\mm{X}\to \mm{T}$

\noindent{ {\bf Игровая задача. }}

  Рассмотрим следующую игровую задачу двух
  игроков.
 Сначала один игрок выбирает множество
 $\Phi \in \mm{V},$ после чего другой игрок
 получает право выбора произвольного
 $x\in \Phi.$

 Пусть есть некоторый критерий качества  $c:C(\mm{T},\mm{X})\to \mm{T}.$
  Пусть один игрок стремится получить
   возможно большее значение
 $c(x),$ %\rav\inf_{t\in\mm{T}}m(t,x(t))$
% для некоторой  непрерывной функции
% качества $m:\mm{T}\times\mm{X}\to \mm{T}$;
 цель другого игрока противоположна.
 Возникает конфликт
  $$ \Phi \in \mm{V} \Uparrow\ \ \  c(x)\ \ \ \  \Downarrow x\in\Phi,$$
%  Качество их выбора при этом можно оценить, например,
%  при помощи введенной на $\mm{K}$ функции
%  $\displaystyle c_{M}(x)\rav \inf(\{t\in\mm{T}\,|\,(t,x(t))\in M\}\cup\{+\infty\}).$
 значением получившейся при этом игры  можно считать максимин
 $$W\rav\sup_{\Phi\in\mm{V}}\inf_{x\in\Phi} c(x)\in\langle-\infty,+\infty].$$
% Всюду далее нас будет интересовать именно тот случай, когда
% $$C_M=\infty.$$

 Сначала построим  такое расширение $\td{\mm{V}}\supset\mm{V}$ этой задачи,
 чтобы
 $$W=\max_{\Phi\in\td{\mm{V}}}\min_{x\in\Phi} c(x),$$
 то есть оба игрока имели оптимальные обобщенные программные управления.

\noindent{{\bf Расширение исходной игровой задачи.}}

 % что-то на тему сепарабельности.

%  Легко проверить,
% что в качестве такой метрики можно взять, например,
%% такой метрикой может быть  например, правилом
% $$\rho(f,g)=\sum_{n\in\mm{N}}\frac{1}{2^n}\frac{||f|_{[0,n]}-g|_{[0,n]}||_{C([0,n],\mm{X})}}{1+||f|_{[0,n]}-g|_{[0,n]}||_{C([0,n],\mm{X})}}. $$
%%  Поскольку $\mm{T}$
%% используется метризуемость $X,$ в принципе свести к [0,n] можно было
%%  и без метризуемости, но тогда пропадет и метрика Хаусдорфа, а это не  есть удобно
 Все элементы $exp \mm{K},$ в том числе и $\mm{K},$ являются компактными
 подмножествами $C(\mm{T},\mm{X})$ (см. \cite[3.4.21]{en}) в силу замкнутости этих множеств,
 а также
 равностепенной ограниченности и равномерной непрерывности множеств вида $\mm{K}_t.$
 Таким образом $\mm{V}\subset (comp)(\mm{K}).$
 Но всякое множество непустых подкомпактов компакта, при оснащении его
 топологией Виеториса, само является предкомпактом \cite[%Теорема
 IV.5.1]{phillvv}.
  Обозначим  через $\td{\mm{V}}$ замыкание множества~$\mm{V},$
  как подмножества $(comp)(\mm{K}),$ в топологии Виеториса,
 тогда $\td{\mm{V}}\in(comp)(comp)(\mm{K}).$

 Теперь определим собственно   расширение исходной игры.  В этой, расширенной игре
 цели игроков не меняются, но
  первый игрок выбирает $\Phi$ уже из $\td{\mm{V}},$ а второй \ --- \
 некоторый элемент $x\in\Phi.$ %, цели же их остались прежними.
 Возникает кофликт
  $$ \Phi \in \td{\mm{V}} \Uparrow\ \ \  c(x)\ \ \ \  \Downarrow x\in\Phi,$$
%  Качество их выбора при этом можно оценить, например,
%  при помощи введенной на $\mm{K}$ функции
%  $\displaystyle c_{M}(x)\rav \inf(\{t\in\mm{T}\,|\,(t,x(t))\in M\}\cup\{+\infty\}).$
 значением  расширенной игры  можно считать максимин
 $$\td{C}=\sup_{\Phi\in\td{\mm{V}}}\inf_{x\in\Phi} c(x).$$

\noindent{{\bf Существование максимина.}}

\begin{teo} \label{teopr}
 %{\bf Теорема }~
 {\it
    Пусть $X$ --- метрическое пространство,
   задан компакт траекторий $\mm{K}\in(comp)(C(\mm{T},\mm{X}))$
  и
  непустой набор пучков $\pust\neq\mm{V}\subset(comp)(\mm{K}).$

 Тогда для всякой непрерывной функции
 качества $c:C(\mm{T},\mm{X})\to \mm{T}$
$$ \max_{\Phi\in\td{\mm{V}}}\min_{x\in\Phi} c(x)=
\sup_{\Phi\in{\mm{V}}}\inf_{x\in\Phi} c(x).$$}
\end{teo}
%Имеет место
%$\displaystyle
 %  C=\td{C}=\max_{\Phi\in\td{\mm{V}}}\min_{x\in\Phi} c(x).$}

\doc

%  Заметим, что функция $c$ непрерывна, действительно
%  функция $m$ непрерывна, следовательно непрерывно для всякого $n\in\mm{N}$ и
%   отображение $x\mapsto c_n(x)=\inf_{t\in [0,n]} m(t,x(t)),$
%   кроме того $c(x)=\lim_{n\to\infty} c_n(x)\geqslant 0,$
%   откуда и следует непрерывность
%   отображения $c.$
%   (в принципе если $\mm{R}$ рассматривать как обратный спектр отрезков $[-n,n],$
%   то для непрерывности ограниченность снизу бы и не понадобилась, но такая топология без введения
%   точек со значением $\infty$ для дальнейшего бесполезна,
%   а если ввести, то такой новый случай отображением $e^{m}$ сводится к
%   расматриваемому)

  Поскольку отображение $c$ непрерывно снизу, то %по принципу компактности
 % Вейерштрасса-Лебега-Бэра
   оно достигает на всяком компакте своего минимального значения, в частности
  отображение
  $A:(comp)(\mm{K})\to\mm{T},$ введенное для всякого $\Phi\in(comp)(\mm{K})$
  по правилу $A(\Phi)=\inf_{x\in\Phi} c(x),$
   совпадает с отображением
  $\Phi\mapsto \min_{x\in\Phi} c(x),$ более того, это отображение
  непрерывно
 % .
 % и также полунепрерывно снизу
  зависит от $\Phi\in\td{U}\in(comp)(comp)(\mm{K}).$ Но тогда
%  во-вторых,
%  непрерывно зависит от $\Phi\in\td{U}\in(comp)(comp)(\mm{K}).$
  на компакте $\td{U}$ оно достигает в некоторой своей точке
   максимального значения, таким образом показано, что
  $$\td{C}=\max_{\Phi\in\td{\mm{V}}} A(\Phi)=\max_{\Phi\in\td{\mm{V}}}\min_{x\in\Phi} c(x).$$

  Покажем, что $\td{C}={C}.$
%  Из вложения $\mm{V}\subset\td{\mm{V}}$ сразу следует $\td{C}\geqslant C.$
%  Покажем $\td{C}\leqslant C.$
  Действительно, $cl \mm{V}=\td{\mm{V}},$
  в частности $\mm{V}$ всюду плотно в $\td{\mm{V}},$ тогда и образ
  $A(\mm{V})$ всюду плотен в $cl A(\td{\mm{V}}),$
  отсюда совпадают верхние грани этих множеств, то есть
  $C=\sup_{\Phi\in\mm{V}} A(\Phi)$ совпадает с
  $\td{C}=\sup_{\Phi\in\td{\mm{V}}} A(\Phi)$
  что и требовалось доказать.\bo

\noindent{ {\bf Характеризация $\td{\mm{V}}.$}}

   Поскольку $\mm{T}$ --- локально компактное хаусдорфово
   пространство, то по
    \cite[3.4.11]{en}
   пространство $\mm{X}^\mm{T}$ с компактно-открытой топологией
   гомеоморфно пределу обратного спектра $\{\mm{X}^{K},\pi^{K'}_{K''},(comp)(\mm{T})\}$
    пространств $\mm{X}^K$ с компактно-открытой топологией
   (Здесь $(comp)(T)$ направлено отношением $\subset$ и
   $\pi^{K'}_{K''}(f)=f|_{K''}$ для всех $f\in \mm{X}^{K'}$).
   Далее, в силу того, что  семейство
   $\{[0,n]\,|\,n\in\mm{N}\}$ конфинально в $(comp)(\mm{T}),$
    из \cite[2.5.11]{en}  $C(\mm{T},\mm{X})$ гомеоморфно пределу
   обратной последовательности
   $\{\mm{X}^{[0,n]},\pi^{[0,n]}_{[0,m]},\mm{N}\}$
    пространств $C([0,n],\mm{X})$
   с топологией равномерной сходимости.
    По \cite[III.1.30]{phillvv}, \cite[2.5.6]{en}
   для  $\mm{K},$ как замкнутого подпространства предела обратного
 спектра,  сужение всех объектов на $\mm{K}$ дает обратный спектр,
 предел которого есть $\mm{K},$ отсюда
   $\mm{K}$ гомеоморфно пределу
   обратного спектра $\mm{S}\rav\{\mm{K}_n,\pi^{[0,n]}_{[0,m]}|_{\mm{K}_n},\mm{N}\}$
   пространств $\mm{K}_n= \{x|_{[0,n]}\,|\,x\in \mm{K} \}$
   в топологии равномерной сходимости, то
   есть
   $$\mm{K}\cong\lim_{\leftarrow}\mm{S}=\lim_{\leftarrow}\{\mm{K}_n,\pi^{[0,n]}_{[0,m]}|_{\mm{K}_n},\mm{N}\}.$$
.

% (Согласно \cite[3.12.26(a)]{en}
% поскольку  $\mm{K},$ не только компакт, но и $T_1$ - пространство, имеет место нормальность
%  $(comp)(\mm{K}).$),
 На
 $(comp)(\mm{K})$
  можно ввести
  метрику Хаусдорфа $H.$ Порожденная ей на $(comp)(\mm{K})$
 топология совпадает  (по \cite[Теорема IV.7.3]{phillvv}) с топологией Виеториса.
 Аналогично можно ввести на
 $(comp)(\mm{K}_n)$
  метрику Хаусдорфа $H_n,$ порожденная ею топология также будет топологией Виеториса.
  % по формулам
  %$$H(\Phi',\Phi'')\rav \max_{x\in\Phi'}\min_{{y\in\Phi''}} ||x-y||_{C(\mm{T},\mm{X})}+
  %\max_{y\in\Phi''}\min_{{x\in\Phi'}} ||x-y||_{C(\mm{T},\mm{X})},$$
  %$$H_n(\Phi'_n,\Phi''_n)\rav \max_{x\in\Phi'_n}\min_{{y\in\Phi''_n}} ||x-y||_{C(\mm{T},\mm{X})}+
  %\max_{y\in\Phi''_n}\min_{{x\in\Phi'_n}} ||x-y||_{C([0,n],\mm{X})}.$$
  %(такая метрика )
  Но в силу \cite[3.12.26(f)]{en},% в энгелькинге достаточно хаусдорфовости пространства нонет слова <<функтор"
  \cite[VII.1.13]{phillvv} функтор $(comp)$ непрерывен
  (то есть функторы $(comp),\lim_{\leftarrow}$ перестановочны)
  в категории компактов, то есть
  $$(comp)(\mm{K})\cong(comp)(\lim_{\leftarrow}\mm{S})=
  \lim_{\leftarrow}\{(comp)(\mm{K}_n),\cdot_n|_{\mm{K}_m},\mm{N}\}.$$
  Отсюда для всякого  множества $F\subset(comp)(\mm{K})$
 согласно \cite[2.5.6]{en} выполнено
 $cl F\cong\lim_{\leftarrow}\{cl (F_n),\cdot_n|_{cl (F_m)},\mm{N}\},$
 в частности, для  $\mm{V}$ имеем
   $\td{\mm{V}}\cong\lim_{\leftarrow}\{cl (\mm{V}_n),\cdot_n|_{cl (F_m)},\mm{N}\},$
   то есть
\begin{gather}
\label{mm}
 \td{\mm{V}}=\{\Phi\in exp\,\mm{K}\,|\,\forall n\in\mm{N}\ \Phi_n\in cl_n(\mm{V}_n)\},
\end{gather}
 где для каждого $n\in\mm{N}$ $ cl_n (\mm{V}_n)$ --- замыкание $\mm{V}_n$
 в компакте $(comp)(\mm{K}_n)$
 (замыкание $\{\Phi_n\,|\,\Phi\in\mm{V}\}$ в метрике
 Хаусдорфа, построенной на всевозможных компактах из $\mm{K}_n$).

% {\bf Замечание} Пусть для всякого $n\in\mm{N}$ найдутся компакт ${Z}^{(n)}$
%  и непрерывное отображение $h^{(n)}$ компакта $Z^{(n)}$
%  на $\mm{V}_n.$
%  , что
%  $h(Z)=\mm{V}$
%  и для всякого $n\in\mm{N}$
%  отображение $z\mapsto(h(z))_n$ непрерывно.
%  Тогда $\mm{V}=\td{\mm{V}}.$

% Действительно, как непрерывный образ компакта, тогда компактно и
% $\mm{V}_n$ для всякого $n\in\mm{N},$ но в этом случае по \rref{mm}
% ${\mm{V}}$ совпадает со своим замыканием
% $\td{\mm{V}},$ что и требовалось.

 \begin{zam} Пусть для всякого $n\in\mm{N}$ найдутся такие компакт ${Z}^{(n)}$
  и полунепрерывное сверху мультиотображение $h^{(n)}$ компакта $Z^{(n)}$
  в подмножества $\mm{K}_n,$ что $h^{(n)}({Z}^{(n)})=\mm{V}_n.$
%  , что
%  $h(Z)=\mm{V}$
%  и для всякого $n\in\mm{N}$
%  отображение $z\mapsto(h(z))_n$ непрерывно.
  Тогда $\mm{V}=\td{\mm{V}}.$
\end{zam}

  Действительно, в силу
  \cite[1.2.35]{Miskis} теперь  для всех $n\in\mm{N}$ компактны
 множества $\mm{V}_n,$ но в этом случае по \rref{mm}
 ${\mm{V}}$ совпадает со своим замыканием
 $\td{\mm{V}},$ что и требовалось.

\noindent{{\bf  Абстрактная задача сближения-уклонения.}}

   Пусть задано некоторое замкнутое в $\mm{T}\times\mm{X}$  целевое множество $M.$
   Пусть игрок, выбирающий траекторию, стремится
 обеспечить $(t,x(t))\in M$ для возможно меньшего $t\in\mm{T},$
 цель выбирающего пучок траекторий противоположна.

   Фактически это означает, что имеет место конфликт
   $$ \Phi \in \mm{V} \Uparrow\ \ \  c_{M}(x)\ \ \ \  \Downarrow x\in\Phi$$
   для функции качества $c:C(T,X)\to \mm{R},$ определенной  по правилу:
   $$\displaystyle c_{M}(x)\rav \inf(\{t\in\mm{T}\,|\,(t,x(t))\in M\}\cup\{+\infty\}).$$
   легко видеть, что эта функция не является  непрерывной в топологии $C(\mm{T},\mm{X}).$

 Введем
  $$
   C_M\rav\sup_{\Phi\in\mm{V}}\inf_{x\in\Phi} c_{{M}}(x).$$

 Если $C_M<\infty,$ то
 имеет место сближение за конечный промежуток времени
 (для некоторого $T\in\mm{T},$ $T\geqslant C_M$ существует такая программа
$\Phi\in\mm{V},$
 что для всех $x\in\Phi$ при некотором $t\in[0,T]$
 выполнено $(t,x(t))\in M.$)

 В противном случае $C_M=\infty,$ то есть
%  Если же  значение   не является точным,  то
% для любых $\Phi\in\mm{V},$ $x\in\Phi$ имеет место
%  $c_{M}(x)<\infty,$ но при этом
 для всех $n\in\mm{N}$
  при некотором $\Phi\in\mm{V}$ для всякой траектории $x\in\Phi$
  выполнено $c_{M}(x)>n.$

 В монографии \cite[c.\,261]{ks} приведен пример такой конфликтно-управляемой
 системы, в которой  (при соответствующей формализации) имеет место
 $C_M=\infty,$ однако нет такого $\Phi\in\mm{V},$ при котором
 для всякой траектории $x\in\Phi$
  выполнено $c_{M}(x)=\infty.$

 Будем говорить, что  {\it значение $C_M=\infty$  является точным}.
 если существует такое $\Phi\in\mm{V},$ что  имеет место уклонение на
 всё $\mm{R},$ то есть
 для всех $x\in\Phi$ при любом $t\in\mm{T}$
 выполнено $(t,x(t))\not\in M,$
  естественно при этом $C_M$ действительно равно $\infty.$

 Понятно, что точное значение $C_M=\infty$ предпочтительней для первого игрока,
 поскольку  гарантирует уклонение сразу на всем временном
 промежутке $\mm{T},$ тогда как в противном случае он может себе
 гарантировать уклонение
 лишь на любом заданном конечном промежутке.

 Рассмотрим теперь также возмущенный вариант задачи, рассмотрим
 некоторое замкнутое в $\mm{T}\times\mm{X}$ множество $\td{M},$
 внутренность которого содержит $M.$
 Тогда аналогично можно рассмотреть
  $$\displaystyle
  c_{\td{M}}(x)\rav \inf\Big(\big\{t\in\mm{T}\,\big|\,(t,x(t))\in
  \td{M}\big\}\cup\big\{+\infty\big\}\Big),\ \
   C_{\td{M}}\rav\sup_{\Phi\in\mm{V}}\inf_{x\in\Phi} c_{\td{M}}(x).$$

 Будем говорить, что {\it значение $C_M=\infty$  является
 грубым (по $M$)}, если для некоторого замкнутого в $\mm{T}\times\mm{X}$
 множества $\td{M},$ $M\subset int \td{M}$ выполнено
 $C_{\td{M}}=\infty.$

 Понятно также, что грубое  значение $C_M=\infty$ предпочтительней для
 первого игрока,
 поскольку  гарантирует уклонение от множества $M$  пусть и с маленьким, но зазором
 (не обязательно равномерно  по $t$ отделенным от нуля).

 Как показывает все тот же пример из \cite{ks}, существует такая
 задача уклонения-наведения, для которой
 $C_M=\infty,$ однако свойство грубости не имеет место.

\begin{teo} \label{teopr1}
%{\bf Теорема}
{ \it Пусть $X$ --- метрическое пространство,
   задан компакт траекторий $\mm{K}\in(comp)(C(\mm{T},\mm{X}))$
  и
  непустой набор пучков $\pust\neq\mm{V}\subset(comp)(\mm{K}).$

Значение $C_M=\infty$ грубо в исходной задаче
  тогда и только тогда, когда значение $\td{C}_M=\infty$ точно
  в расширенной задаче.}
\end{teo}
\doc

  Пусть $C_M=\infty$ точно
  в расширенной задаче. Тогда для некоторого $\td\Phi\in\td{\mm{V}}$
  выполнено $(t,x(t))\not\in M$ для всех $t\in\mm{T},x\in\td\Phi.$
%  В частности $(0,x_0)\not\in M.$

   Введем множество $Y$ всех $x\in\mm{K},$ для которых выполнено $(t,x(t))\in M$
   хотя бы при каком-то $t\in\mm{T}.$ Легко видеть, что его дополнение открыто, а тогда
   само множество $Y$
   замкнуто.

%
%      Определим $Y\rav\{(t,x(t))\,|\,t\in\mm{T},\,x\in\td\Phi\}.$
 %   В силу точности оно не пересекается с $M.$
%    Заметим, что
%    для любого $n\in\mm{N}$ $\td\Phi_n\in cl(\mm{V}_n)\subset(comp)(\mm{K}_n).$
%    Но тогда и замкнуто множество
%    $$Y_n\rav Y\cap
%    ([0,n]\times\mm{X})=\{(t,x(t))\,|\,t\in\mm{T},\,x\in\td\Phi_n\}.$$
%    Поскольку всякое пересечение $Y$ с любым компактным в
%    $\mm{T}\times\mm{X}$ множеством    теперь замкнуто,
%    а пространство $\mm{T}\times\mm{X},$ как и всякое
%    метрическое пространство  является хаусдорфовым секвенциальным, %\cite[с.371]{en},
%    а следовательно $k$-пространством \cite[3.3.20]{en},
%     то замкнуто и само $Y$ согласно \cite[3.3.18]{en}.

  Поскольку $Y$ и $\td{\Phi}$ замкнуты в метрическом (а, следовательно, совершенно нормальном \cite[4.1.13]{en})
  пространстве $C(\mm{T},\mm{X}),$ то по критерию совершенной нормальности
  Веденисова
  \cite[1.5.19]{en}
  найдется такая непрерывная функция
  $c:C(\mm{T},\mm{X})\to[0,1],$ что $c^{-1}(0)=Y,$ $c^{-1}(1)=\td\Phi.$
%  Назначим функции $z:\mm{T}\times\mm{X}\to\mm{T}$ и
%  $m:\mm{T}\times\mm{T}\to \mm{T}$ правилами
%  $z(t,x)=a(x), m(t,r)=r$ для всех $x\in\mm{X},t,r\in\mm{T}.$

 %  Заметим, что функция $c$ непрерывна, действительно
%   Поскольку  функция $m$ непрерывна, то непрерывно для всякого $n\in\mm{N}$ и
%   отображение $x\mapsto c_n(x)\rav\inf_{t\in [0,n]} m(t,x(t)),$
%   теперь из ограниченности снизу отображения
%   $x\mapsto c(x)\rav\lim_{n\to\infty} c_n(x)\geqslant 0,$
%   следует  и его непрерывность.
%   (в принципе если $\mm{R}$ рассматривать как обратный спектр отрезков $[-n,n],$
%   то для непрерывности ограниченность снизу бы и не понадобилась, но такая топология без введения
%   точек со значением $\infty$ для дальнейшего бесполезна,
%   а если ввести, то такой новый случай отображением $e^{m}$ сводится к
%   расматриваемому)

  Рассмотрим игру с функцией качества $c,$ тогда
 по уже показанному предложению
 $$%\sup_{\Phi\in{\mm{V}}}\inf_{x\in\Phi}\inf_{t\in\mm{T}}a(t,x(t))
  %C=\mm{C}=
  \sup_{\Phi\in{\mm{V}}}\inf_{x\in\Phi}c(x)=\max_{\Phi\in\td{\mm{V}}}\min_{x\in\Phi}c(x)\geqslant\min_{x\in\td\Phi}c(x)=1,$$
 отсюда найдется такое $\Phi'\in\mm{V},$ что
  $\min_{x\in\Phi'}c(x)>1/2.$

  Теперь рассмотрим
  $${M}'\rav\{(t,x(t))\in\mm{T}\times\mm{X}\,|\, x\in c^{-1}([0,1/2\rangle),t\in\mm{T}\}.$$
  Как прообраз, открыто $c^{-1}([0,1/2\rangle),$ из непрерывности графика, определенной
  на компакте функции, теперь для всех $n\in\mm{N}$  открыты
  множества ${M}'\cap
    ([0,n]\times\mm{X}),$ тогда это верно и для $M'.$
  Введем  множество $\td{M}\rav cl {M}'$.
  Из определения $c$ следует, что $M\subset{M}'\subset int \td{M}.$
  Но, по только что показанному, найдется $\Phi'\in\mm{V},$
  для которого  выполнено
  $c(x)>1/2,$ отсюда
  $(t,x(t))\not\in \td{M}$
  при всех $x\in\Phi',$ $t\in\mm{T}.$
  Таким образом выполнено $C_{\td{M}}=\infty$ и
  значение   $C_{M}=\infty$  грубо в исходной задаче.

  Пусть исходная задача груба для некоторого $\td{M}.$ Рассмотрим
  множество $D=cl \big(\mm{T}\times(\mm{X}\setminus\td{M})\big),$ по условию
  это  множество не пересекается с $M.$
%  По лемме Урысона
% \cite[I.2.18]{phillvv}
%  найдется такая непрерывная функция
%  $a:\mm{T}\times\mm{X}\to[0,1],$ что $a^{-1}(0)=M,$ $a^{-1}(1)=N.$
%  Назначим функции $z:\mm{T}\times\mm{X}\to\mm{T}$ и
%% $z(t,x)=z(t,x), m(t,r)=r$ для всех $x\in\mm{X},t,r\in\mm{T}.$

  Теперь для всякого $n\in\mm{N}$
  найдется такое $\Phi^n\in\mm{V},$ что
для всех $t\in[0,n],x\in\Phi^n$ выполнено  $(t,x(t))\not\in\td{M},$
 то есть $(t,x(t))\in D.$
 Поскольку $\Phi^n$ из метрического компакта $\mm{V},$ можно
 (при необходимости перейдя к подпоследовательности) считать,
 что имеется предел $\td{\Phi}=\lim_{n\to\infty}\Phi^n\in\mm{V}.$
 Тогда для всех $t\in[0,n],x\in\td\Phi$ выполнено
 $(t,x(t))\in cl D=D,$ в частности $(t,x(t))\not\in M.$
 Тогда $\td{C}_M$ точно. \bo

 %{\bf Пример 1.}
 % Определим для всех $p\in (0,1]$ функцию $y_p:\mm{R}_{\geqslant 0}\to\mm{R}$
 %по правилу $y_p(x)=(x-p)^2/p^2$ для всех
% $x\in\mm{R}_{\geqslant 0}.$ Определим $U\rav\{\{y_p\}\,|\,p\in(0,1]\}.$
% Пусть $M\rav\mm{R}\times{0}.$ Легко видеть, что $c_{\td{M}}(x)=\infty,$
% это значение грубо, но неточно. Расширим игру, тогда
% $\td{U}=u\cup\{y_0\},$ где  $y_0:\mm{R}_{\geqslant0}\to\mm{R}$ --- поточечный предел $y_p$ при $p\to
% 0,$ то есть $y_0(x)= 1$ для всех
% $x\in\mm{R}_{\geqslant 0}.$
% В расширенной задаче $y_0$ гарантирует уклонение на весь промежуток
% времени.

\begin{pr} \label{prpr1}% {\bf Пример 1.}
  Пусть $$\Phi=\Big\{x(t)\in C\big(\mm{T},\{r\}\big)\,\big|\,
  r\in\langle-1,0\rangle\cup\langle0,1\rangle\Big\},\
  \ \mm{V}=\{\Phi\},\  M\rav\mm{T}\times\{0\}.$$
  Легко видеть, что в такой <<игре>> имеет место точность, но нет
  грубости. Но это и неудивительно, ведь множества $\mm{V}_t$ не замкнуты.
  \end{pr}

\begin{pr} \label{prpr2}% {\bf Пример 2.}
  Определим для всех $p\in (0,1]$ функцию $y_p:\mm{R}_{\geqslant 0}\to\mm{R}$
 по правилу $y_p(x)=(x-p)^2/p^2$ для всех
 $x\in\mm{R}_{\geqslant 0}.$ Определим $U\rav\{\{y_p\}\,|\,p\in\langle0,1]\}.$
 Пусть $M\rav\mm{R}\times{0}.$ Легко видеть, что $c_{\td{M}}(x)=\infty,$
 это значение грубо, но неточно. Расширим игру, тогда
 $\td{U}=U\cup\{y_0\},$ где  $y_0:\mm{R}_{\geqslant 0}\to\mm{R}$ --- поточечный предел $y_p$ при $p\to
 0,$ то есть $y_0(x)= 1$ для всех
 $x\in\mm{R}_{\geqslant 0}.$
 В расширенной задаче $y_0$ гарантирует уклонение на весь промежуток
 времени, то есть расширенная задача точна; понятно также, что малое
  изменение $M$ не может изменить $C_M,$ равного $+\infty.$
\end{pr}

%{\bf Следствие 1.}
\begin{sle} \label{slepr1}
{\it Если $\td{\mm{V}}={\mm{V}},$ то значение $V_M=\infty$ грубо
  тогда и только тогда, когда  точно.}
\end{sle}

\begin{sle} \label{slepr2}
%{\bf Следствие 2.}
{\it Если ${\mm{V}}$ конфинально в  $(\td{\mm{V}},\supseteq),$
  то значение $V_M=\infty$ грубо
  тогда и только тогда, когда  точно.}
\end{sle}

  Для доказательства достаточно заметить, что в условиях
  конфинальности если расширенная задача точна, то такова и исходная.
  Действительно, если расширенная задача точна, то найдется
  пучок  $\td\Phi\in\td{\mm{V}},$ гарантирующий $c_M(x)=\infty$ для всякой
  траектории $x\in\td\Phi.$
  По свойству конфинальности  в ${\mm{V}}$  имеется пучок $\Phi,$ для
  которого
  $\td\Phi\supseteq\Phi.$ Но тогда тем более
  $c_M(x)=\infty$ для всякой
  траектории $x\in\Phi,$ то есть исходная задача точна.

 \begin{zam} Пока в этом подразделе целевое множество $M$ было замкнутым в $\mm{X}$
  множеством, однако можно взять произвольное замкнутое в $C(\mm{T},\mm{X})$
  множество $Y,$ определить для всякого $x\in C(\mm{T},\mm{X})$ момент выхода из $Y$
  $$c_{\mm{Y}}(x)=\inf (\{t\in\mm{T}\,|\,x|_{[0,t]}\in Y_t\})$$
  и, понимая <<грубость>> как нечувствительность максимина при
  уменьшении $Y$ до некоторого  содержащегося во внутренности $Y$ замкнутого
  множества $Y'$, показать эквивалентность  <<точности>> исходной и <<грубости>> расширенной
  игры. Это можно доказать непосредственно, повторив соответствующие рассуждения теоремы,
   а можно, переформулировав показанную теорему
  для пространства
  $C(\mm{T},\mm{X})$ в качестве нового $\mm{X}$.
\end{zam}

\noindent{{\bf Игровая задача сближения-уклонения с дискриминацией
первого игрока}}

 Рассмотрим конфликтно-управляемую систему
\begin{gather}\label{sys_}
 \dot{x}=f(t,x,u,v),\ x(0)=x_0, u\in P,\ q\in Q,
\end{gather}
 функционирующую   в
 конечномерном  пространстве $\mm{X}$  на бесконечном промежутке времени $\mm{T}.$
  В \rref{sys_} $P$  и $Q$ ---  компакты в
 конечномерных пространствах.
 Функция $f:
 \mm{T} \times \mm{X} \times P \times Q  \mapsto
 \mm{X}$ предполагается:
 1) непрерывной; 2) локально липшицевой по второй (фазовой)
 переменной, то есть для любого $A\in(comp)(\mm{T}\times \mm{X})$ найдется такая
 суммируемая функция $L_A\in B(\mm{T},\mm{R}),$ что
 $$||f(t,x_1,u,v)-f(t,x_2,u,v)||\leqslant L_A(t)||x_1-x_2||,\ \ \forall
 (t,x_1),(t,x_2)\in A, u\in P, v\in Q; $$
 3) удовлетворяющей условию подлинейного роста, то есть
 при некоторой локально суммируемой функции $a\in B(\mm{T})$ имеет место
 $$  ||f(t,x,u,v)||_{\mm{X}}\leqslant a(t)(1+||x||_{\mm{X}})\ \ \ \forall (t,x,u,v)\in \mm{T}\times\mm{X}\times P\times Q.$$

Рассмотрим теперь дифференциальное включение
\begin{gather}\label{vkl}
\dot{x}\in \overline{co} f(t,x,P,Q),\  x(0)=x_0.
\end{gather}
% $$\dot{x}\in co f(t,x,P,Q),\  x(0)=x_0.$$
 Заметим, что в силу условия подлинейного роста
  всякое его локальное решение продолжимо вправо
 \cite{tovst},
 обозначим через $\mm{K}$ множество всех его продолжимых вправо
 решений. Тогда для всякого $t\in\mm{T}$ множество $\mm{K}_t$
 действительно будет компактным (как и полагается множеству, обозначенному через
  $\mm{K},$ в первом разделе).
  Для всякого семейства функций $S,$ определенных на $\mm{T},$
  условимся через $S_n$ обозначать всевозможные их сужения на $[0,n].$

\noindent{{\bf Иллюстрация равенства $\td{\mm{V}}=\mm{V}.$}}

 Пусть {$\Pi(\mm{T}\times Q)$}
  --- множество всех регулярных
  борелевских мер $\eta$ на {$\mm{T}\times Q$}
 % (на {$\mm{T}\times P\times Q$}),
  имеющих лебеговскую проекцию на
  элементы $(comp)(\mm{T}),$
  то есть для всякого $I\in(comp)(\mm{T})$
  $\eta(I\times Q)=\lambda (I).$
  Аналогично введем множество  {$\Pi(\mm{T}\times P \times Q)$}
  всех  регулярных
  борелевских мер $\eta$
  на {$\mm{T}\times P\times Q$}, для которых
  $\eta(I\times P \times Q)=\lambda (I)$ при
  любом $I\in(comp)(\mm{T}).$

  Каждой  $\eta\in {\Pi(\mm{T}\times P \times Q)}$ можно сопоставить
  $\ph(\eta)\in C(\mm{T},\mm{X})$ как
  решение уравнения
   $$x(t)=x_0+\int_{[0,t]\times P\times Q} f(\tau,x(\tau),u,v)\,d\eta,\ \ \forall t\in\mm{T}.$$
  Легко видеть, что это решение из $\mm{K}.$

  Для всех ${\nu\in\Pi(\mm{T}\times Q)}$ введем
  множество ${\Pi_{{\nu}}(\mm{T}\times P\times Q)}$
   всех регулярных
  борелевских  на $T\times P\times Q$
 мер $\eta\in{\Pi(\mm{T}\times P \times Q)},$\  согласованных с мерой $\nu,$ то есть
  таких, что
     для любых борелевских $I\subset\mm{T},A\subset P$ выполнено
     $\eta(I\times P\times A)=\nu(I\times A).$

 После выбора  вторым игроком меры {$\nu\in\Pi(\mm{T}\times
 Q)$} первый игрок выбирает меру
  {$\eta\in\Pi_{{\nu}}(\mm{T}\times P\times Q)$}, то есть
   траекторию ${x=\ph(\eta)\in\ph(\Pi_{{\nu}}(\mm{T}\times P\times
  Q))}.$

  Определим
   $$\displaystyle \mm{V}\rav \{{\ph(\Pi}_{{\nu}}{(\mm{T}\times P\times
  Q))} \,|\, \ {\nu\in\Pi(\mm{T}\times Q)}\}.$$

 Как и на всем $\mm{T}$, для  отрезка $[0,n]$
 (при всяком $n\in\mm{N}$) можно также ввести
 множество {$\Pi([0,n]\times Q)$} всех регулярных
  борелевских мер $\eta$ на {$[0,n]\times Q$},
 % (на {$\mm{T}\times P\times Q$}),
  имеющих лебеговскую проекцию на
  элементы $(comp)([0,n])$;
   множество  {$\Pi([0,n]\times P \times Q)$}
  всех  регулярных  борелевских мер $\eta$
  на {$[0,n]\times P\times Q$}, а кроме того, для всякой
 ${\nu\in\Pi([0,n]\times Q)}$
  множество ${\Pi_{{\nu}}([0,n]\times P\times Q)}$
   всех регулярных
  борелевских  на $[0,n]\times P\times Q$
 мер $\eta\in{\Pi([0,n]\times P \times Q)},$\  согласованных с мерой $\nu.$
  Каждой  $\eta\in {\Pi([0,n]\times P \times Q)}$ аналогично можно сопоставить
  $\ph(\eta)\in C([0,n],\mm{X})$.

  Заметим, что
  символ $\cdot|_n$ пропускается <<внутрь>> этих множеств, то есть:
 $$\Pi(\mm{T}\times Q)_n=\Pi([0,n]\times Q),\
  \Pi(\mm{T}\times P\times Q)_n=\Pi([0,n]\times P\times Q)$$
  $$   \Pi_{{\nu}}(\mm{T}\times P\times Q)_n=
   \Pi_{{\nu}|_{[0,n]}}([0,n]\times P\times Q),\ \ \ \forall\nu\in\Pi(\mm{T}\times Q).$$
  Поскольку это  касается также траекторий и образов отображения $\ph$, то
 выполнено
  $$\mm{V}_n=\{{\ph(\Pi}_{{\nu}}{([0,n]\times P\times
  Q))} \,|\, \ {\nu\in\Pi([0,n]\times Q)}\}.$$

  Как показано в \cite{Real,dep}  отображение
 $$\nu\in\Pi([0,n]\times
  Q)\mapsto \ph(\Pi_{{\nu}}([0,n]\times P\times
  Q))\subset \mm{K}_n$$
  имеет компактные значения и,
  как многозначное отображение, непрерывно по Хаусдорфу,
 теперь
  множество $\mm{V}_n\subset(comp)(\mm{K}_n)$ ---
  непрерывный образ компакта  $\Pi([0,n]\times Q),$ следовательно тоже компакт.
  Поскольку это верно для всякого $n\in\mm{N},$
  то из замечания 2 следует, что $\mm{V}=\td{\mm{V}}.$
  Теперь из следствия 1 получаем

\begin{utv} \label{utvpr}
{\it %{\bf Следствие 3}{
  В паре <<программное смешанное  управление\ /\ программное обобщенное  управление>>
     значение $C_M=\infty$ точно тогда и только тогда, когда грубо.
}
\end{utv}

 Отметим, что то же имеет место в паре <<позиционная стратегия\ /\ позиционная
 стратегия>>: значение $C_M=\infty$ точно тогда и только тогда, когда грубо.
 Доказательство этого фактически содержится в доказательстве \cite[Теорема
 60.1]{ks}. При этом, впрочем, не строится никаких топологических расширений
 множества позиционных стратегий, вместо этого достаточно рассмотреть неограниченно
 возрастающую последовательность стабильных мостов и отображение,
 сопоставляющее каждому мосту стратегию для игры $\mm{V}_n$

\vspace{3ex}

\small

\makeatletter \@addtoreset{equation}{section}
\@addtoreset{footnote}{section}
\renewcommand{\section}{\@startsection{section}{1}{0pt}{1.3ex
plus 1ex minus 1ex}{1.3ex plus .1ex}{}}

{ %\scriptsize

\renewcommand{\refname}{{\rm\centerline{СПИСОК ЛИТЕРАТУРЫ}}}

%, the projective limit
%  compact-open topology.
%Pursuit-evasion game
%topology of Vietoris an extension of the  original game the game be
%played under <<joint control\ /\ relaxed control".
%  The value $W_M=\infty$ is
% rough   iff the value $W_M=\infty$  is
%  precise.
%Максиминные программы Программное управление
\titleeng

\end{document}